\documentclass{rhumj_new}

\usepackage{algorithm}
\usepackage{algorithmic}
\usepackage{bm}
\usepackage{tikz}
\usetikzlibrary{arrows.meta}

\title[Keep the beat going]{Keep the beat going: Automatic drum transcription with momentum}

\author[A.L.~Foster]{Alisha L. Foster}
\affiliation{University of California San Diego}
\email{a1foster@ucsd.edu}
\author{Project advisor Robert J. Webber}
\affiliation{University of California San Diego}
\email{rwebber@ucsd.edu}


\begin{document}
\abstract{How can we process a piece of recorded music to detect and visualize the onset of each instrument?
A simple, interpretable approach is based on partially fixed nonnegative matrix factorization (NMF). 
Yet despite the method's simplicity, partially fixed NMF is challenging to apply because the associated optimization problem is high-dimensional and non-convex.
This paper explores two optimization approaches that preserve the nonnegative structure, including a multiplicative update rule and projected gradient descent with momentum.
These techniques are derived from the previous literature, but they have not been fully developed for partially fixed NMF before now.
Results indicate that projected gradient descent with momentum leads to the higher accuracy among the two methods, and it satisfies stronger local convergence guarantees.}

\section{Introduction}

Drums are the rhythmic backbone of a song, contributing to our perception and enjoyment of music.
A typical drum kit consists of many parts--- such as a snare drum, bass drum, and multiple cymbals--- which are played simultaneously to create a drum groove. 
This paper considers the computational inverse problem: given a recorded musical piece, how can we extract the drum signals that created it?

Automatic drum transcription (ADT) is the act of extracting and notating drum grooves from recorded music.
In ADT, the data are transformed into \emph{spectrograms}, time-frequency representations calculated with the short-time Fourier transform.
The spectrograms are analyzed via ADT systems to precisely detect the timing of the drum hits.
Yet a challenge for ADT systems is that drums are inharmonic instruments, meaning they do not produce pure musical tones like a piano or guitar. 
Moreover, each drum or cymbal sounds different based on physical aspects such as its size, material, and tension.
This makes drum transcription more difficult than the transcription of a piano, for example, where each A4 note can be fully characterized by its fundamental frequency of 440Hz.
A typical ADT system needs to be adapted to the individual drum kit.

\begin{figure}[t]
\centering

\begin{minipage}{0.46\linewidth}
    \centering
    \includegraphics[width=\linewidth]{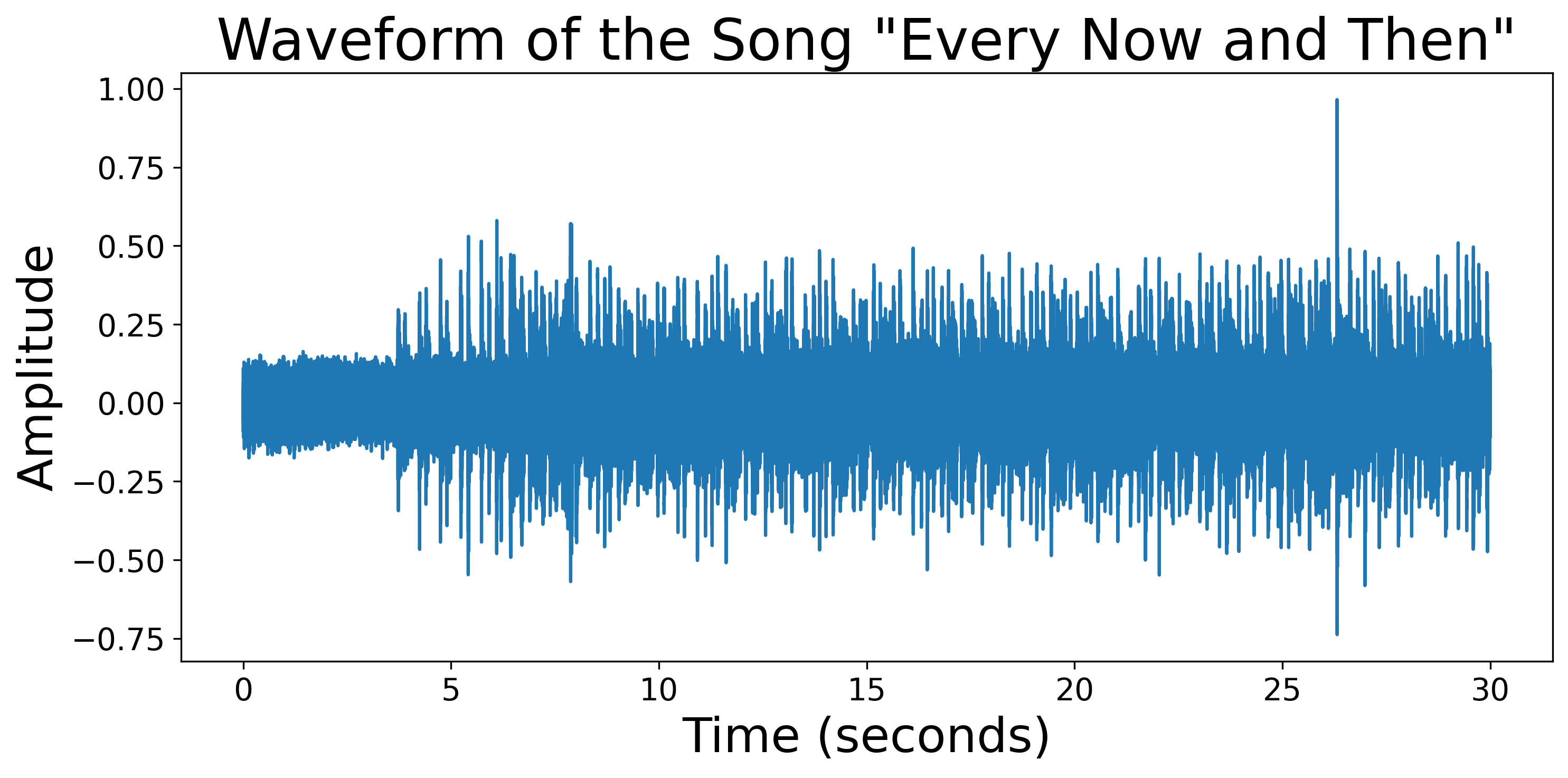}
\end{minipage}
\begin{minipage}{0.03\linewidth}
    \centering
    \begin{tikzpicture}
        \draw[line width=4pt,->,>={Stealth[length=10pt,width=10pt]}]
            (0,0) -- (0.75,0);
    \end{tikzpicture}
\end{minipage}
\hfill
\begin{minipage}{0.46\linewidth}
    \centering
    \includegraphics[width=\linewidth]{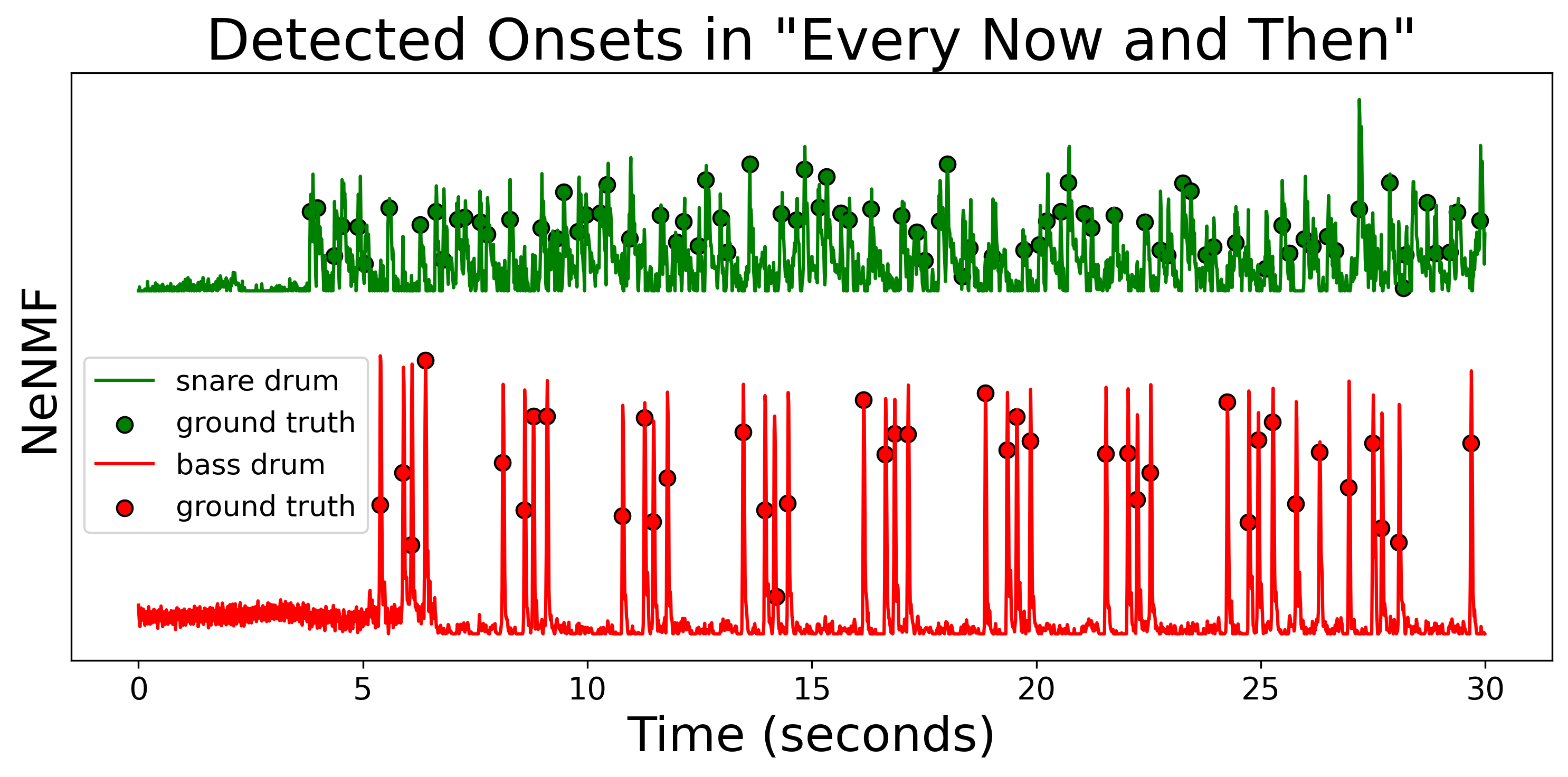}
\end{minipage}
\caption{Automatic drum transcription (ADT) system applied to the first 30 seconds of the song ``Every Now and Then'' to visualize the onsets of snare and bass drums.}
\label{fig:placeholder}
\end{figure}

Ultimately, an ADT system could be useful to drummers as a tool for notating their grooves or learning beats from recorded musical pieces.
The ideal ADT system would be fully interpretable, 
identifying the frequency representations of each part of the drum kit.
This work constructs an interpretable ADT method based on partially fixed nonnegative matrix factorization (NMF), and it identifies an efficient numerical method for the partially fixed NMF optimization.
See \cref{fig:placeholder} for an illustration of the interpretable ADT system constructed in this work.
The optimization approach is based on established methods \cite{lee2000algorithm,GuanNeNMF}, but the theory and implementation of these methods need to be modified for the partially fixed NMF problem.
The rest of this introductory section describes the partially fixed NMF problem and the optimization approach in more detail.

\subsection{Constructing an interpretable ADT model}
To give a brief background, NMF approximates a nonnegative-valued data matrix $\bm{V}$ as the product of two nonnegative-valued matrices, called a dictionary matrix $\bm{W}$ and an activation matrix $\bm{H}$.
\begin{equation*}
    \bm{V} \approx \bm{W} \bm{H}.
\end{equation*}
In the context of automatic music transcription, 
$\bm{V}$ is the $m \times n$ magnitude spectrogram matrix that has $m$ frequency bins and $n$ time frames.
$\bm{W}$ is an $m \times r$ dictionary matrix that contains a basis of frequencies for the song, and
$\bm{H}$ is an $r \times n$ activation matrix that contains the temporal onsets of each frequency component.
In this context, it is natural that $\bm{V}$ is approximated by the product of $\bm{W}$ and $\bm{H}$ matrices that both have nonnegative entries.

Paatero and Tapper advocated for NMF in 1994 and called it ``positive matrix factorization'' \cite{paatero1994positive}.
Later, Lee and Seung popularized the method in 2000 \cite{lee2000algorithm}.
NMF has many applications in music analysis \cite{musicanalysis}, image representation \cite{imagerep}, recommender systems \cite{recsys}, and genomics \cite{genomics}. 
\Cref{sec:litreview} reviews existing works that apply NMF to ADT.

\begin{figure}[H]
    \centering
    \includegraphics[width=\linewidth]{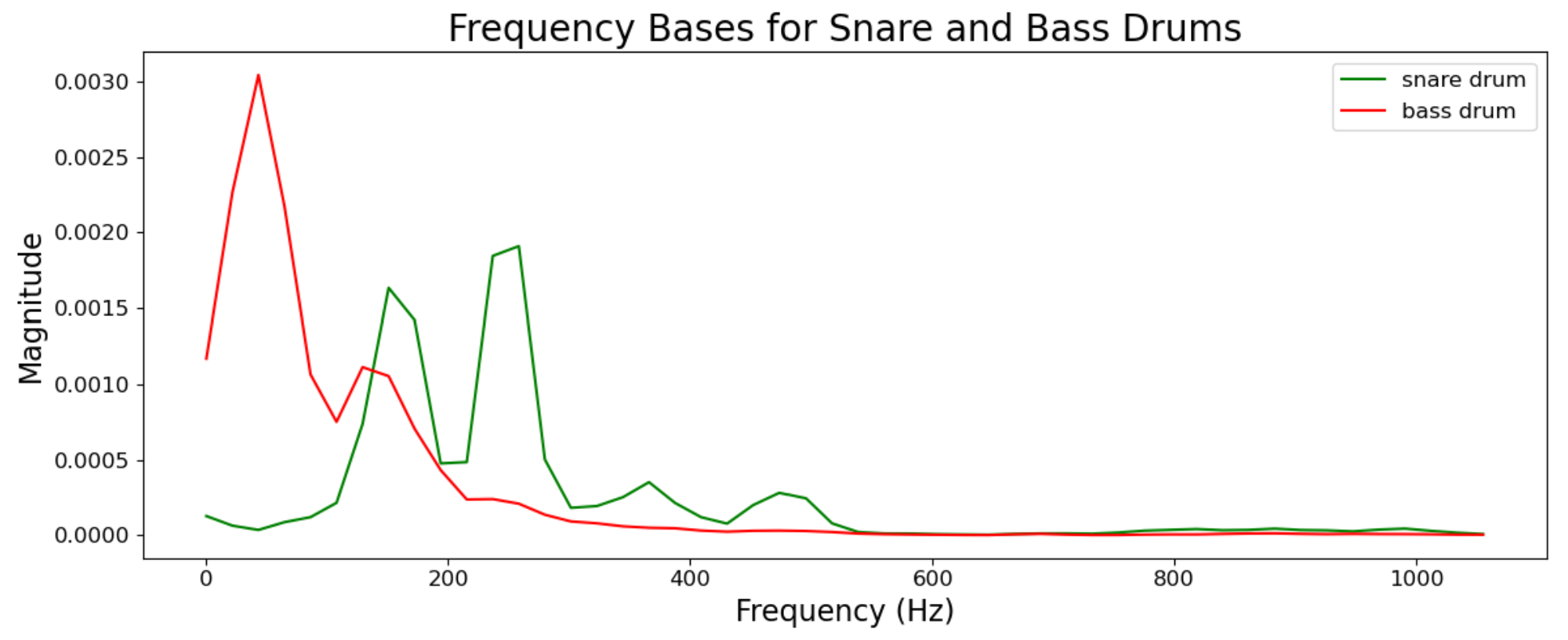}
    \caption{Plot of the magnitudes of the frequencies for the basis components in $\bm{W}_D$ for the snare and bass drums for ``Every Now and Then."}
    \label{fig:new_fig}
\end{figure}

As the main limitation, the conventional NMF model is underspecified. 
Multiple columns of $\bm{W}$ can represent the same instrument or multiple instruments can be represented in the same column.
Therefore, our work will investigate a more highly specified model called partially fixed NMF \cite{YooEtAl,WuLerch}, where the dictionary and activation matrices are factored as
\begin{equation*}
\bm{W} = \begin{bmatrix}
    \bm{W}_D \quad \bm{W}_H
\end{bmatrix} \quad \text{and} \quad
\bm{H} = \begin{bmatrix}
    \bm{H}_D \\
    \bm{H}_H
\end{bmatrix}.
\end{equation*}
Here $\bm{W}_D$ is an $m \times r_D$ matrix,
where the rank $r_D$ is the number of drum parts to be transcribed, and each column contains one drum component that we want to learn, e.g., snare drum, bass drum, and hi-hat cymbal.
In practice, a drummer could calibrate $\bm{W}_D$ by recording some drum hits to create a frequency basis for each part of a drum kit.
For example, \cref{fig:new_fig} illustrates the columns of $\bm{W}_D$ corresponding to snare and bass drums in green and red lines.
This set of columns stays fixed for the duration of the optimization process.
$\bm{H}_D$ is an $r_D \times n$ matrix that contains the temporal onsets of the known drum instruments needed for transcription.
For example, \cref{fig:placeholder} illustrates the rows of $\bm{H}_D$ corresponding to snare and bass drums.
Last, $\bm{W}_H$ contains the frequencies of the harmonic instruments in the song, and the rank $r_H$ is tuned to balance accuracy with interpretability.
Similarly, $\bm{H}_H$ is the $r_H \times n$ activation matrix for the harmonic instruments of the song. 
$\bm{W}_H\bm{H}_H$ essentially absorbs all non-drum noise from the spectrogram. 

\subsection{Addressing the optimization challenge}

It is challenging to optimize the objective function in the partially fixed NMF model, which is given by
\begin{equation*}
    F(\bm{H}_D, \bm{W}_H, \bm{H}_H) = \frac{1}{2} \lVert \bm{V} - (\bm{W}_D\bm{H}_D + \bm{W}_H\bm{H}_H) \rVert_{\rm F}^2.
\end{equation*}
The objective function is high-dimensional and non-convex.
Nonetheless, the NMF community has identified an optimization approach using alternating minimization that works well in practice.
Freeze two matrices in the factorization and optimize the other matrix, which becomes a convex optimization problem.
Then trade the roles of the matrices.
Repeat as many times as desired.
\begin{align*}
    \bm{H}_D^{t+1} &= \operatornamewithlimits{argmin}_{\bm{H}_D \geq 0} F(\bm{H}_D, \bm{W}_H^t, \bm{H}_H^t) \\
    \bm{W}_H^{t+1} &= \operatornamewithlimits{argmin}_{\bm{W}_H \geq 0} F(\bm{H}_D^{t+1}, \bm{W}_H, \bm{H}_H^t) \\
    \bm{H}_H^{t+1} &= \operatornamewithlimits{argmin}_{\bm{H}_H \geq 0} F(\bm{H}_D^{t+1}, \bm{W}_H^{t+1}, \bm{H}_H).
\end{align*}
There is no guarantee the approach will converge to a global minimizer, but the square Frobenius norm error can be made to decrease at each step.

The optimization part of this paper focuses on a subproblem in the alternating minimization:
how do we optimize one matrix when the others are fixed?
There are at least two natural strategies that preserve the nonnegative structure: (1) a multiplicative update rule based on gradient descent \cite{lee2000algorithm} and (2) a projected gradient descent approach with Nesterov momentum \cite{GuanNeNMF}.
We need to execute these methods carefully to make sure the iterates remain nonnegative, and we need to carefully summarize the relevant theory to support each approach.
This work will adapt the two optimization strategies to the problem of partially fixed NMF, developing theory and algorithms.
The results show that projected gradient descent with Nesterov momentum has favorable convergence guarantees and outperforms the multiplicative update rule in its empirical accuracy with a fixed runtime.

\subsection{Outline}

The paper is organized as follows. 
\Cref{sec:litreview} reviews existing work, \cref{sec:background} describes the partially fixed NMF problem and formulates optimization methods, \cref{sec:drum} applies the methods to automatic drum transcription of the ENST-Drums database \cite{ENST} and a demo of ``Every Now and Then" by the author's band, Threadbare \cite{song3}, and \cref{sec:conclusion} concludes.

\subsection{Notation}

Scalars are in regular typeface: $n, d, L$.
Vectors are in bold lower-case letters: $\bm{x}, \bm{y}$.
Matrices are in bold capital letters: $\bm{M}, \, \bm{N}$.
The $i$th entry of a vector $\bm{x}$ is $\bm{x}(i)$.
The $(i,j)$ entry of a matrix $\bm{M}$ is $\bm{M}(i,j)$, the $i$th row is $\bm{M}(i, \cdot)$, and the $j$th column is $\bm{M}(\cdot, j)$.
The transpose is $\bm{M}^\top$, the Moore-Penrose pseudoinverse is $\bm{M}^+$, the spectral norm is $\lVert \bm{M} \rVert$, and the Frobenius norm is $\lVert \bm{M} \rVert_{\rm F}$.
The space of nonnegative-valued $m \times n$ matrices is denoted $\mathbb{R}^{m \times n}_+$.
The element-wise product of two vectors $\bm{x}, \bm{y}$ is $\bm{x} \odot \bm{y}$.
Similarly, the element-wise product of two matrices $\bm{M}, \, \bm{N}$ is $\bm{M} \odot \bm{N}$, and the element-wise division is $\frac{\bm{M}}{\bm{N}}$.

\section{Literature Review}
\label{sec:litreview}

This section reviews ADT models from the literature. 
Drawbacks to some existing models are the large amount of data required for optimization, the slow optimization speed, and the lack of human interpretability. 

\subsection{A variety of ADT models}

ADT has been performed with several neural network-based models.
For example, Jacques and Roebel \cite{CNN} trained three convolutional neural networks using a data set of over 400 musical pieces --- one each for detecting the hi-hat, snare drum, and bass drum --- and achieved a high degree of accuracy.
As a drawback, however, neural networks are not the most interpretable model choice.
Therefore, many researchers have applied NMF for more interpretable transcriptions of drum music \cite{SmaragdisConvNMF,YooEtAl,WuLerch} or polyphonic music \cite{Polyphonic, SemiSupCNMF, OST}.
NMF takes advantage of the inherent nonnegativity and low-rank properties of magnitude spectrograms, and a possible advantage of NMF models is that they require little training data.

One popular version of the NMF model, called ``convolutive NMF'' (CNMF), uses a convolution operator to account for local dependencies between nearby time frames.
Smaragdis \cite{SmaragdisConvNMF} used the CNMF model to transcribe monophonic drum data.
Later, H. Wu et al. \cite{SemiSupCNMF} used CNMF 
to transcribe piano music.
The CNMF method is not perfect, however.
A drawback noted by H. Wu et al. is the long computation time needed for transcriptions, although recent work has pursued more robust optimization algorithms for CNMF \cite{GillisCNMF}.

Any NMF model can be used with a variety of loss functions. The Frobenius norm loss function is common \cite{lee2000algorithm, GuanNeNMF, PGD4NMF}, and two other popular loss functions in audio source separation are the Kullback-Leibler divergence \cite{GillisBSS, lee2000algorithm} and Itakura-Saito divergence \cite{GillisBSS, musicanalysis}.
The latter two loss functions are logarithmically scaled, mirroring human audio perception \cite{GillisBSS}. 
Here we use the simplest, most common Frobenius norm loss function, but exploring alternative loss functions is an interesting topic for future work.

\subsection{Partially fixed NMF}

As a limitation, traditional NMF models are underspecified, so they might assign multiple columns of the dictionary matrix to the same instrument or assign one column to multiple instruments.
In order to transcribe sheet music in this context, a drummer would need to train an auxiliary model identifying which instruments produce which onsets.
This work focuses on a more highly specified NMF model called ``partially fixed NMF'' that was proposed by Yoo et al. \cite{YooEtAl} and later extended by C.-W. Wu and Lerch \cite{WuLerch}.
The major benefit of the partially fixed NMF is interpretability.
In partially fixed NMF, each drum instrument is assigned to a single column of the dictionary matrix, and the drum onsets are described by the corresponding rows of the activation matrix.
For example, the first column of the dictionary matrix might represent the snare drum, the second column might represent the bass drum, etc. and then the first row of the activation matrix would contain the snare drum onsets and the second row would contain the bass drum onsets.
This ordering is helpful for transforming the data into human-readable drum notation.

\subsection{Optimizers for NMF}

There are many optimizers for NMF problems, the most classic and prevalent being the multiplicative update rule (MUR) presented by Lee and Seung \cite{lee2000algorithm}.
Projected gradient descent is another first-order method suitable for high-dimensional optimization \cite{PGD4NMF},
and Guan et al. have developed a projected gradient descent method with Nesterov momentum for NMF called ``NeNMF'' \cite{GuanNeNMF}.
Other NMF optimizers include alternating least squares \cite{paatero1994positive}, alternating nonnegative least squares \cite{ANLS}, the alternating direction method of multipliers \cite{ADMM}, hierarchical alternating least squares \cite{GillisHALS}, and proximal block coordinate descent \cite{GillisBlock}.
Here we focus on analyzing, applying, and comparing MUR and NeNMF.
As motivation, MUR is the most classic and widely used optimizer, while NeNMF converges more rapidly than MUR in the experiments of \cite{GuanNeNMF} and \cite[Sec.~8.4]{gillis}.
MUR and NeNMF are highly flexible NMF optimization methods that can be adapted to a range of optimization problems, as emphasized in \cite[Tab.~8.2]{gillis}.
Therefore, it makes sense to develop these methods first when optimizing the partially fixed NMF.
Nonetheless, comparisons between more optimizers would be interesting further work.

\section{Methods and optimization techniques}
\label{sec:background}

This section will describe the MUR and NeNMF optimization methods in detail and derive the best available mathematical convergence guarantees, building on techniques in \cite{lee2000algorithm,GuanNeNMF}.

\subsection{Multiplicative update rule}
Lee and Seung \cite{lee2000algorithm} studied the classic NMF problem
\begin{equation}
\label{NMF}
\min_{\bm{W}, \bm{H}} \frac{1}{2}  \lVert \bm{V} -\bm{WH} \rVert_{\rm F}^2,
\end{equation}
and they proved the following result.
\begin{thm}[Convergence of original MUR algorithm \cite{lee2000algorithm}]
\label{thm:MUR}
  The square Frobenius norm $\lVert \bm{V} -\bm{WH} \rVert_{\rm F}^2$ is nonincreasing under the multiplicative update rule
\begin{equation}
  \label{eq:original_mur}
      \bm{H} \leftarrow \bm{H} \odot \frac{\bm{W}^\top\bm{V}}{\bm{W}^\top\bm{WH}}, \qquad \bm{W} \leftarrow \bm{W} \odot \frac{\bm{VH}^\top}{\bm{WHH}^\top}.
\end{equation}
\end{thm}

In a critique of the MUR method, Lin \cite{MURconvergence} pointed out that terms in the denominator might equal zero, in which case the method is ill-defined.
However, zero denominator terms only arise in pathological settings, and
\cref{thm:MUR} remains valid as long as the denominators stay positive.

In 2010, Yoo et al. \cite{YooEtAl} proposed solving the partially fixed NMF problem
\begin{equation}
\label{PFNMF}
    \min_{\bm{H}_D, \bm{W}_H, \bm{H}_H} \frac{1}{2} \lVert \bm{V} - (\bm{W}_D\bm{H}_D + \bm{W}_H\bm{H}_H) \rVert_{\rm F}^2.
\end{equation}
using a variant of the multiplicative update rule \cref{eq:original_mur}.
Their approach is a little subtle, since we cannot simply replace $\bm{W}$ with $\bm{W}_D$, $\bm{H}$ with $\bm{H}_D$, and $\bm{V}$ with $\bm{V} - \bm{W}_H \bm{H}_H$ in \cref{eq:original_mur} to produce an update
\begin{equation*}
	\bm{H}_D \leftarrow \bm{H}_D \odot \frac{\bm{W}_D^\top[\bm{V} - \bm{W}_H \bm{H}_H]}{\bm{W}_H^\top\bm{W}_H \bm{H}_D}
\end{equation*}
We tried this and ran into numerical problems, since the iterates might take on negative entries!
To preserve nonnegativity, we need to be more careful and apply the approach in \cref{alg:MUR}, derived from the work \cite{YooEtAl}.

As a limitation, the paper \cite{YooEtAl} only outlined a justification of \cite[Thm.~1]{lee2000algorithm} following the original proof strategy, but they never proved the method's validity in detail.
As another issue, the proof strategy in \cite[Thm.~1]{lee2000algorithm} assumes positive-valued entries, but the algorithm should remain valid as long as the denominator terms are positive.
To fill in these gaps, \cref{proof} establishes a full theoretical guarantee for the MUR method that uses a slightly different argument from  \cite[Thm.~1]{lee2000algorithm} based on evaluating the best multiplicative (rather than additive) perturbations.
This argument leads to the following result.

\begin{thm}[Convergence of MUR algorithm for partially fixed NMF]
\label{thm:PFMUR}
    The loss function 
    $F(\bm{H}_D, \bm{W}_H, \bm{H}_H) 
    = \frac{1}{2} \| \bm{V} - (\bm{W}_D\bm{H}_D + \bm{W}_H\bm{H}_H) \|^2_F
    $
    is nonincreasing under the multiplicative update rule
    \begin{equation*}
    \begin{aligned}
        \bm{H}_D &\leftarrow \bm{H}_D \odot \frac{\bm{W}_D^\top\bm{V}}{\bm{W}_D^\top(\bm{W}_D\bm{H}_D + \bm{W}_H\bm{H}_H)}, \\
        \bm{W}_H &\leftarrow \bm{W}_H \odot \frac{\bm{VH}_H^\top}{(\bm{W}_D\bm{H}_D + \bm{W}_H\bm{H}_H)\bm{H}_H^\top}, \\
        \bm{H}_H &\leftarrow \bm{H}_H \odot \frac{\bm{W}_H^\top\bm{V}}{\bm{W}_H^\top(\bm{W}_D\bm{H}_D + \bm{W}_H\bm{H}_H)},
    \end{aligned}
    \end{equation*}
    as long as the denominator terms are positive.
\end{thm}

\Cref{alg:MUR} uses the multiplicative update rule established in \cref{thm:PFMUR}, preserving the nonnegativity of the iterates and the monotonicity of the loss function.

\begin{algorithm}[H]
    \caption{MUR for partially fixed NMF}\label{alg:MUR}
    \begin{algorithmic}[1]
        \STATE Construct $\bm{W}_D \in \mathbb{R}_+^{m \times r_D}$, where $r_D$ is the number of drum components
        \STATE Randomly initialize $\bm{H}_D^1 \in \mathbb{R}_+^{r_D \times n}, \bm{W}_H^1 \in \mathbb{R}^{m \times r_H}_+, \bm{H}_H^1 \in \mathbb{R}^{r_H \times n}_+$
        \FOR{$t=1,...,T$}
            \STATE $\bm{H}_D^{t+1} =  \bm{H}_D \odot \frac{\bm{W}_D^\top\bm{V}}{\bm{W}_D^\top(\bm{W}_D\bm{H}_D^t + \bm{W}_H^t\bm{H}_H^t)}$
            \STATE $\bm{W}_H^{t+1} = \bm{W}_H^t \odot \frac{\bm{V} (\bm{H}_H^t)^\top}{(\bm{W}_D\bm{H}_D^{t+1} + \bm{W}_H^t\bm{H}_H^t){(\bm{H}_H^t)^\top}} $
            \STATE $\bm{H}_H^{t+1} = \bm{H}_H^t \odot \frac{({\bm{W}_H^{t+1}})^\top\bm{V}}{({\bm{W}_H^{t+1}})^\top(\bm{W}_D\bm{H}_D^{t+1} + \bm{W}_H^{t+1}\bm{H}_H^t)}$
        \ENDFOR
        \STATE $\bm{H}_D = \bm{H}_D^T$, $\bm{W}_H = \bm{W}_H^T$, $\bm{H}_H = \bm{H}_H^T$.
        \STATE return $\bm{H}_D, \bm{W}_H, \bm{H}_H$
    \end{algorithmic}
\end{algorithm}

\subsection{Projected gradient descent with momentum}

In 2012, Guan and coauthors \cite{GuanNeNMF} introduced the ``NeNMF'' method for the classic NMF problem \cref{NMF}.
Like MUR, the NeNMF method is based on an alternating minimization over $\bm{W}$ and $\bm{H}$ factors.
As the major innovation, NeNMF solves each subproblem
\begin{equation*}
\min_{\bm{H} \geq 0} \lVert \bm{W} \bm{H} - \bm{V} \rVert_{\rm F}^2
\quad \text{or} \quad
\min_{\bm{W} \geq 0} \lVert \bm{W} \bm{H} - \bm{V} \rVert_{\rm F}^2
\end{equation*}
to high accuracy.
To do this, NeNMF applies several iterations of projected gradient descent with Nesterov momentum, which Guan et al. call the ``optimal gradient method'' (OGM) and which is presented in pseudocode in
\cref{alg:OGMHD}.

\begin{algorithm}[H]
\caption{$\operatorname{OGM}(\bm{W}, \bm{H}_0, \bm{V})$}
\label{alg:OGMHD}
\begin{algorithmic}[1]
    \STATE Initialize $\alpha_1 = 1$, $\bm{Y}_1 = \bm{H}_0$
    \FOR{$k=1, ..., K$}
        \STATE $\bm{H}_k = \operatorname{ReLu} \Bigl( \bm{Y}_k - \frac{\bm{W}^\top [ \bm{W} \bm{Y}_k - \mathbf{V} ]}{\lVert \bm{W} \rVert^2} \Bigr) $
        \STATE $\alpha_{k+1} = \frac{1 + \sqrt{4\alpha_k^2 + 1}}{2}$
        \STATE $\bm{Y}_{k+1} = \bm{H}_k + \frac{\alpha_k - 1}{\alpha_{k+1}} (\bm{H}_k - \bm{H}_{k-1})$
    \ENDFOR
    \STATE $\bm{H} = \bm{H}_K$
    \STATE return $\bm{H}$
\end{algorithmic}
\end{algorithm}

OGM is based on gradient descent, and it projects each iterate onto the cone of nonnegative-valued matrices by an element-wise application of the rectified linear unit $\operatorname{ReLu}(x) = \max\{x, 0\}$.
The method incorporates Nesterov momentum using a step size parameter $\alpha_k$, which promotes an accelerated convergence rate.
Guan et al. established the following guarantee for OGM in \cite[Prop.~2.1]{GuanNeNMF}.
\begin{thm}[Convergence of OGM \cite{GuanNeNMF}]
\label{thm:OGM}
    Let $\bm{W} \in \mathbb{R}^{m \times r}$, $\bm{H}_0 \in \mathbb{R}_+^{r \times n}$, and $\bm{V} \in \mathbb{R}^{m \times n}$ be inputs to the optimal gradient method (\cref{alg:OGMHD}), and fix a matrix
    \begin{equation*}
        \bm{H}_{\star} \in \operatornamewithlimits{argmin}_{\bm{H} \geq \bm{0}} \lVert \bm{W} \bm{H} - \bm{V} \rVert_{\rm F}^2.
    \end{equation*}
    Then, after $K$ iterations, the optimal gradient method
    outputs a matrix $\bm{H} = \bm{H}_K$ that satisfies
    \begin{equation*}
        \lVert \bm{W} \bm{H}_K - \bm{V} \rVert_{\rm F}^2
        \leq \lVert \bm{W} \bm{H}_{\star} - \bm{V} \rVert_{\rm F}^2 + \frac{2\lVert \bm{W} \rVert^2 \lVert \bm{H}_0 - \bm{H}_\star \rVert_{\rm F}^2}{(K+2)^2}.
    \end{equation*}
\end{thm}
\Cref{thm:OGM} ensures that OGM solves each subproblem up to an additive error of $\mathcal{O}\big(\frac{1}{K^2}\big)$ after $K$ iterations.

To our knowledge, NeNMF has not previously been applied to the partially fixed NMF problem \cref{PFNMF}.
However, we can extend NeNMF to the partially fixed NMF problem by using OGM to solve each subproblem in the alternating minimization, as presented in
\cref{alg:NeNMF}.
This new NeNMF extension is straightforward, with none of the conceptual or theoretical difficulties as in MUR.
Because \cref{alg:NeNMF} is using OGM to solve each subproblem, convergence guarantees follow immediately from \cref{thm:OGM}.

\begin{algorithm}[H]
\caption{NeNMF for partially fixed NMF}\label{alg:NeNMF}
\begin{algorithmic}[1]
    \STATE Construct $\bm{W}_D \in \mathbb{R}_+^{m \times r_D}$, where $r_D$ is the number of drum components
    \STATE Randomly initialize $\bm{H}_D^1 \in \mathbb{R}_+^{r_D \times n}, \bm{W}_H^1 \in \mathbb{R}^{m \times r_H}_+, \bm{H}_H^1 \in \mathbb{R}^{r_H \times n}_+$
    \FOR{$t=1,...,T$}
        \STATE $\bm{H}_D^{t+1} = \operatorname{OGM}(\bm{W}_D, \bm{H}_D^t, \bm{V} - \bm{W}_H^t \bm{H}_H^t) $
        \STATE $\bm{W}_H^{t+1} = \operatorname{OGM}((\bm{H}_H^t)^\top, (\bm{W}_H^t)^\top, (\bm{V} - \bm{W}_D \bm{H}_D^{t+1})^\top)^\top$
        \STATE $\bm{H}_H^{t+1} = \operatorname{OGM}(\bm{W}_H^{t+1}, \bm{H}_H^t, \bm{V} - \bm{W}_D \bm{H}_D^{t+1})$
    \ENDFOR
    \STATE $\bm{H}_D = \bm{H}_D^T$, $\bm{W}_H = \bm{W}_H^T$, $\bm{H}_H = \bm{H}_H^T$.
    \STATE return $\bm{H}_D, \bm{W}_H, \bm{H}_H$
\end{algorithmic}
\end{algorithm}

\subsection{Comparison of theoretical guarantees and runtimes}

Obtaining the global minimizer of the partially fixed NMF problem \cref{PFNMF} is challenging since the loss function is high-dimensional and nonconvex.
Nevertheless, we can hope to solve each convex subproblem in the alternating minimization to high accuracy.
NeNMF is guaranteed to solve each subproblem up to an additive error of $\mathcal{O}(1/K^2)$ due to \cref{thm:OGM}.
In contrast, MUR is not guaranteed to solve each subproblem to any level of accuracy, although it does decrease the overall loss due to \cref{thm:PFMUR}.
Therefore, NeNMF has stronger theoretical guarantees.

Updating a single matrix in an iteration of MUR or OGM requires $\mathcal{O}(mnr)$ arithmetic operations where $r = r_D + r_H$ is the total rank of drum components and harmonic components.
Since NeNMF performs $K$ iterations of OGM while MUR performs a single multiplicative update for each matrix, the cost of one iteration of NeNMF is roughly $K$ times higher than the cost of one iteration of MUR.
The results in \cref{sec:drum} account for this factor-of-$K$ difference by running MUR for $T = 100$ iterations versus running NeNMF with $T = 10$ outer iterations containing $K = 10$ inner OGM iterations.
Therefore, the runtime between algorithms is similar (65 seconds for MUR, 68 seconds for NeNMF).

\section{Application to automatic drum transcription} \label{sec:drum}

This section describes numerical experiments applying partially fixed NMF to automatic drum transcription.
\Cref{sec:data} describes the data sets and \cref{sec:results} presents results.

\subsection{Data sets} \label{sec:data}

ENST-Drums \cite{ENST} is a publicly available data set containing annotated audio recordings from three drummers playing individual hits, short phrases, and longer songs in various styles on their own drum kits.
We downloaded the \textit{minus-one} subset of the ENST-Drums data from the repository \cite{sungkyun2023yourmt3}.
This subset contains 28 tracks of length $\approx 70$ s with the drummers playing their drum kits without accompaniment from other instruments.
The sample rate is 16 kHz.

\begin{figure}
    \centering
    \includegraphics[width=0.8\linewidth]{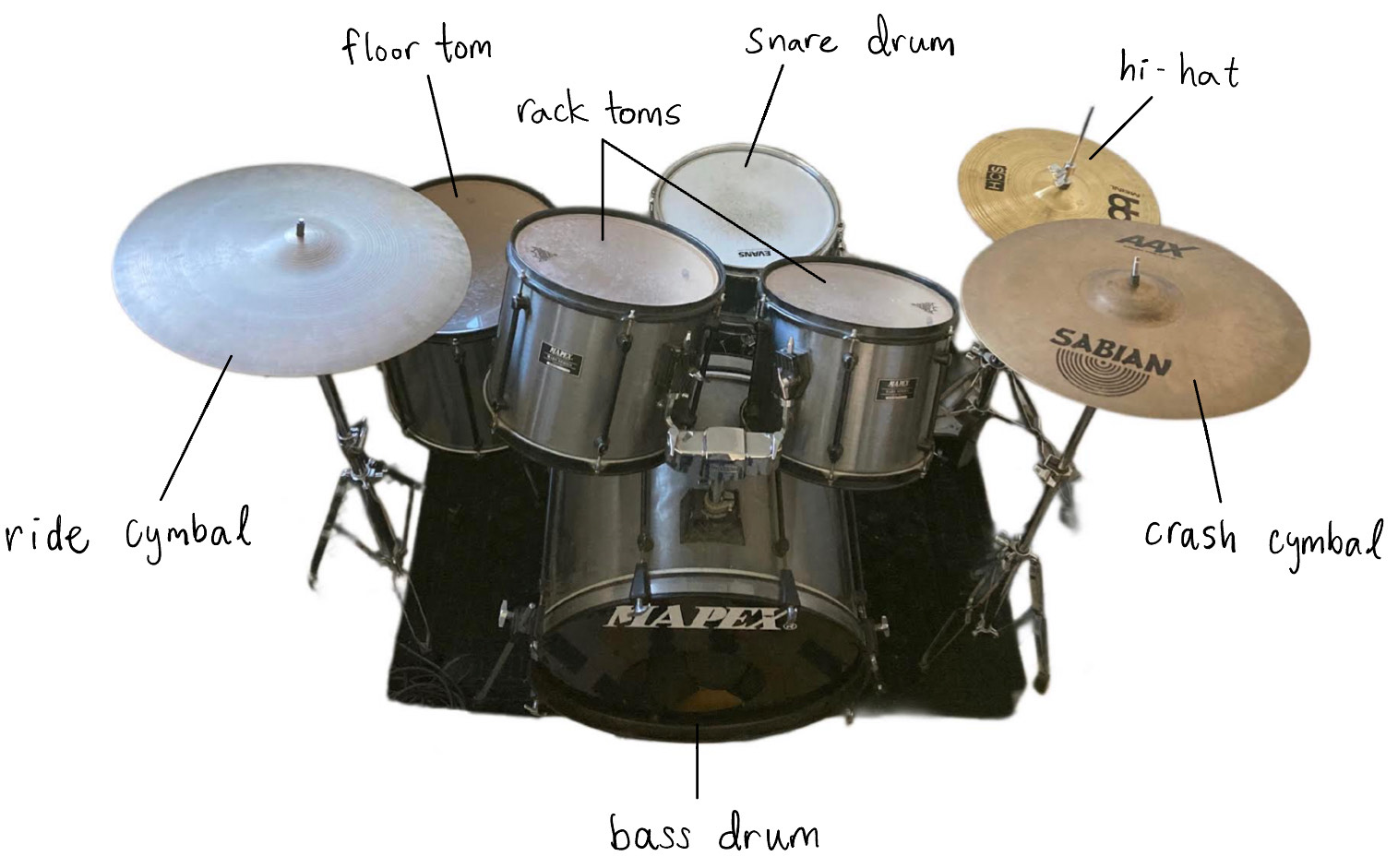}
    \caption{Photograph of the author's drum kit with instruments annotated by hand.
    The snare and bass drums were included in the automatic drum transcription for the song ``Every Now and Then'', as described in \cref{sec:drum}.}
    \label{fig:drum}
\end{figure}

``Every Now and Then" is a recorded song from the author's band, Threadbare \cite{song3}.
It consists of two guitars, one bass guitar, vocals, and a drum kit.
It is available for streaming at \url{https://open.spotify.com/track/3FoWRPLO9Z0tQiVY9jLrth?si=1e252b92cde04a8e}. 
The sample rate is 44.1 kHz. 
We manually annotated each onset of the snare drum or bass drum as a ground-truth reference.
See \cref{fig:drum} for a labeled photo of the complete drum kit, including the two annotated components.

As a preprocessing procedure, we converted ENST-Drums and ``Every Now and Then'' into magnitude spectrograms using the short-time Fourier transform with a window length 2048 and hop size 512. 
This transformation resulted in a time resolution of 32 ms for ENST-Drums and 12 ms for ``Every Now and Then". The spectrogram size for ``Every Now and Then" was 1025 frequency bins by 13950 time frames. The spectrograms of the tracks in ENST-Drums also included 1025 frequency bins, with varying lengths around 2200 time frames.

Last, we constructed a dictionary element for each annotated drum component by extracting the magnitude spectrogram of an individual hit and averaging over the time axis.
We constructed dictionary elements for the hi-hat, snare drum, and bass drum for each of the three drum kits in ENST-Drums.
We constructed dictionary elements for the snare drum and bass drum for ``Every Now and Then''.

\subsection{Results} \label{sec:results}

We applied the MUR and NeNMF algorithms to solve the partially fixed NMF problem \cref{PFNMF} for ENST-Drums and ``Every Now and Then'' data using a harmonic rank $r_H=5$.
In these tests, we independently initialized the entries of $\bm{H}_D, \bm{W}_H,$ and $\bm{H}_H$ as $\operatorname{Unif}(0,1)$ random variables. 
Then we applied either MUR for 100 iterations or  NeNMF for 10 outer iterations using 10 inner OGM iterations.
These experiments led to similar runtimes for the two algorithms. For ``Every Now and Then," the MUR runtime was about 65 seconds, and the NeNMF runtime was about 68 seconds.

\begin{figure}[H]
    \centering
    \includegraphics[width=0.89\linewidth]{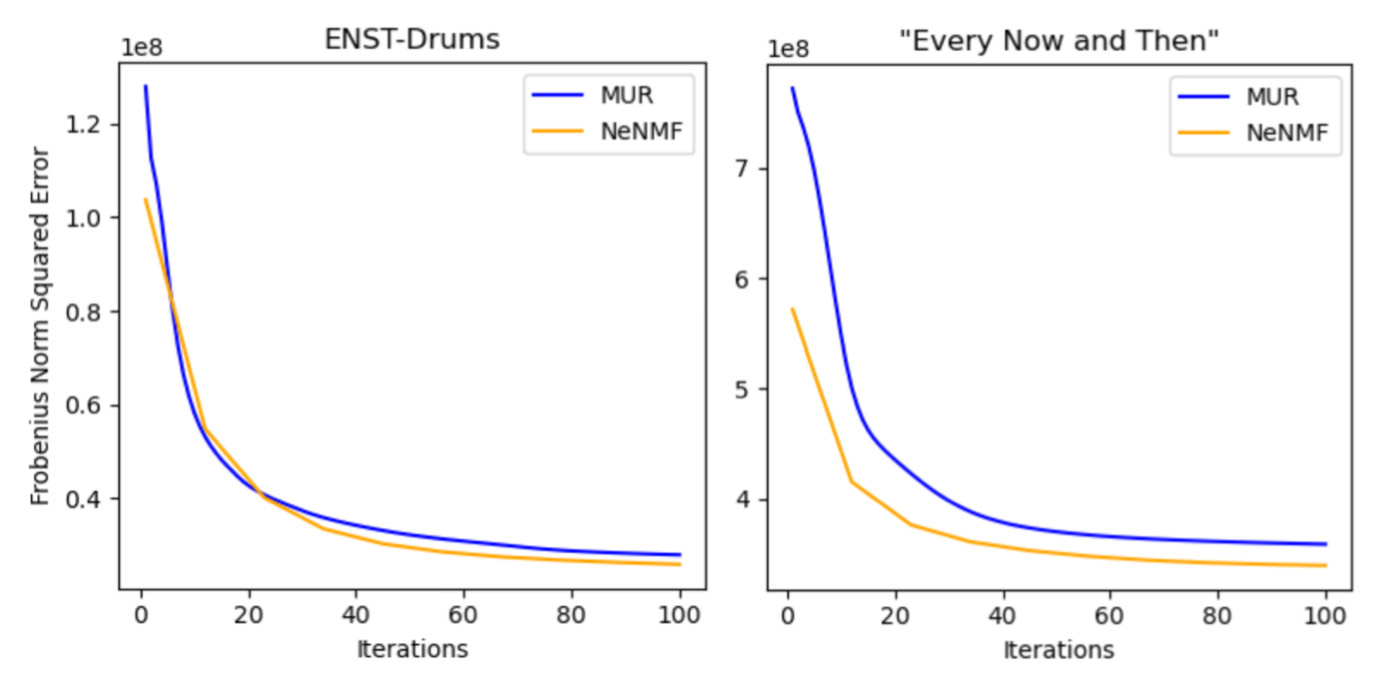}
    \caption{Square Frobenius norm error as a function of runtime for MUR and NeNMF.
    The horizontal axis shows the number of iterations for MUR and the number of inner iterations ($10\times$ the number of outer iterations) for NeNMF. Note that the number of inner OGM iterations per outer iteration for NeNMF is fixed at 10.
    The error for ENST-Drums (left) is averaged across 28 tracks, while the error for ``Every Now and Then" (right) is for a single track.}
    \label{fig:convergence}
    \vspace{0.2cm}
\end{figure}

\Cref{fig:convergence} compares the square Frobenius norm errors achieved using MUR and NeNMF.
The errors for ENST-Drums (left panel) are averaged across the 28 tracks, while the errors for ``Every Now and Then'' (right panel) are for a single track.
To promote a fair comparison, the horizontal axis records the iteration parameter $t$ in MUR (\cref{alg:MUR}) and it records $10\times$ the iteration parameter $t$ in NeNMF (\cref{alg:NeNMF}).
Note that the number of inner OGM iterations for NeNMF is fixed at 10 for these experiments.
The results show that NeNMF leads to a comparable or smaller square Frobenius norm error for either data set and any fixed runtime.

\begin{figure}
    \centering
    \includegraphics[width=\linewidth]{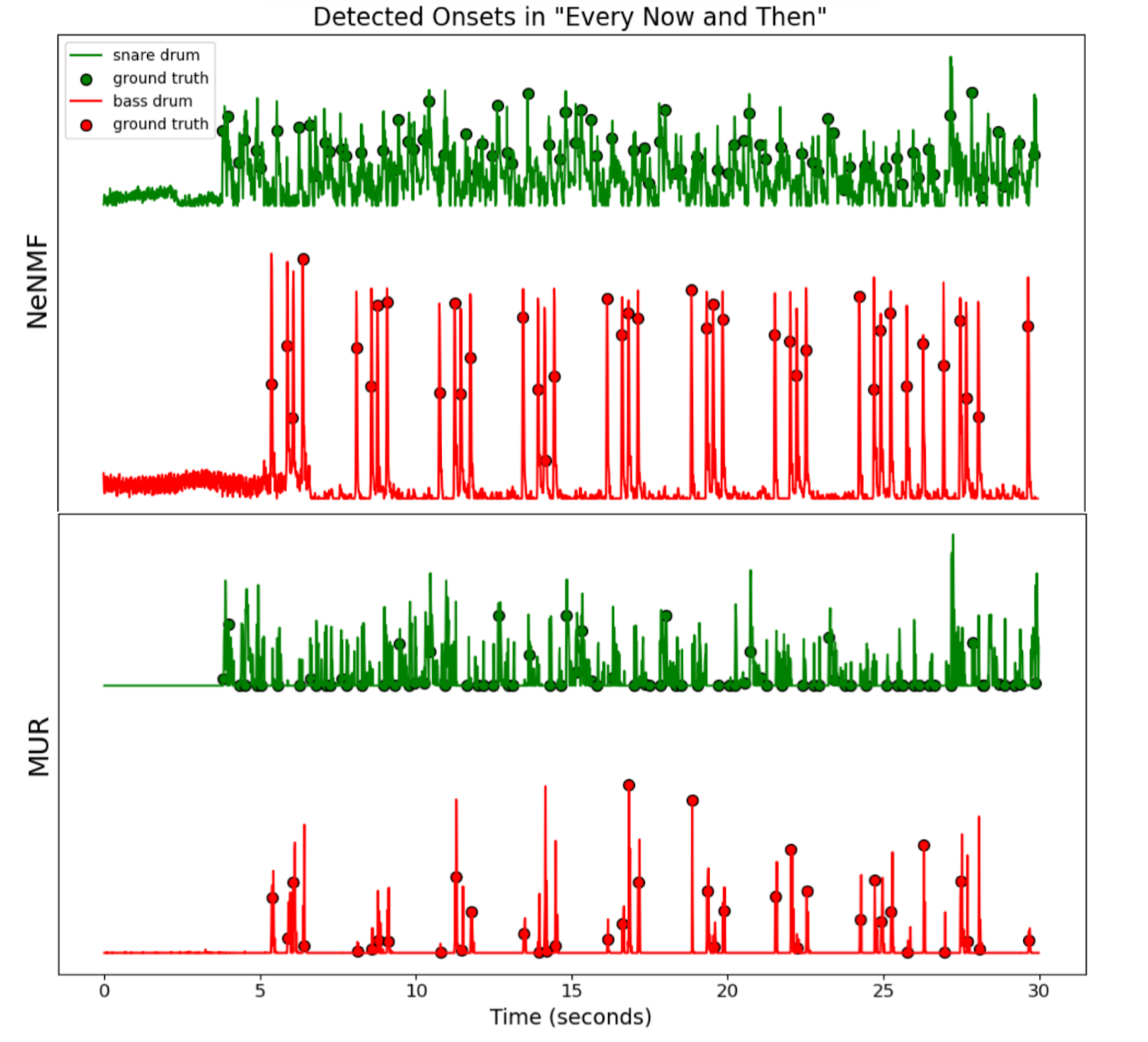}
    \caption{Entry values of drum activation matrix $\bm{H}_D$ plotted along time, optimized with NeNMF (top) and MUR (bottom) for the first 30 seconds of ``Every Now and Then." The markers represent ground-truth annotations, and large entry values of $\bm{H}_D$ represent detected onsets.}
    \label{fig:song3annotation}
\end{figure}

Next, we evaluated the accuracy of the drum hits identified using MUR and NeNMF.
Both optimization algorithms produce drum activation matrices consisting of nonnegative entries $\bm{H}_D(i,j)$, where a large entry represents a detected onset of drum instrument $i$ at time frame $j$. \Cref{fig:song3annotation} visually compares the accuracy of the onsets detected by MUR and NeNMF for the first 30 seconds of ``Every Now and Then." The plots show that $\bm{H}_D$ optimized with NeNMF more accurately captures the drum onsets.

As a complication, there are various ways to determine which entries of $\bm{H}_D$ are large enough to be regarded as drum hits or ``onsets''.
One standard method used in automatic music transcription \cite{OST, daniel:inria-00452615}
identifies $p_j$ of the largest entries in each column $\bm{H}_D(\cdot, j)$, where $p_j$ is the ground-truth number of drum hits.
Although the ground-truth number of drum hits is not likely to be available in real-world use cases, this approach allows for straightforward comparisons without the need to tune evaluation hyperparameters.
After identifying the drum onsets using the ground-truth number of onsets, we counted the true positives (TP), false positives (FP), and false negatives (FN) in each row of the onset matrix $\bm{H}_D$.
Then, we measured accuracy using the $F$-score:
\begin{equation*}
    F = \frac{2{\rm TP}}{2{\rm TP} + {\rm FP} + {\rm FN}}.
\end{equation*}
Each $F$-score lies in the range $F \in [0, 1]$, with higher $F$-scores indicating greater accuracies.

\Cref{table:Eval-1} reports the F-scores calculated using the ground-truth number of offsets for each data set and optimization method, averaged over all the drum components and tracks for each data set.
The results show that NeNMF achieves higher accuracy than MUR for both ENST-Drums and ``Every Now and Then''.
NeNMF achieves an especially high accuracy of 97.5\% for ``Every Now and Then'', because the track contains just two annotated drum components and the onset matrix correctly indicates which component is larger at nearly all the onset times.

\begin{table}[t]
\centering
\label{table:Eval-1}
\begin{tabular}{|c|c|c|}
\hline
 & MUR & NeNMF \\
\hline
\textit{ENST-Drums} & 0.599 & \textbf{0.620} \\
\hline
\textit{``Every Now and Then"} & 0.854 & \textbf{0.975}  \\
\hline
\end{tabular}
\vspace{0.2cm}
\caption{\label{table:Eval-1}F-scores evaluated using ground-truth number of offsets.
Higher F-scores are marked in bold.}
\end{table}

A different way to identify drum hits is a signal-adaptive median threshold \cite{WuLerch, Lerch}.
This method evaluates the median in a sliding window over each row $\bm{H}_D(i, \cdot)$, and it adds a small offset to determine the threshold for a drum hit.
For the experiments, we set the window length to 0.1 s, and we set the offset coefficients to 0.05 for the hi-hat, 0.1 for the snare drum, and 0.15 for the bass drum.
We considered a detected onset correct if it is within 50 ms of a ground-truth onset. 
Note that changing the window length and the offset coefficients can lead to large changes in the evaluation scores, and different tracks have better scores for different values.

\begin{table}[t]
\centering
\begin{minipage}{0.45\textwidth}
    \begin{tabular}{|c|c|c|}
    \hline
     & MUR & NeNMF \\
    \hline
    Hi-hat & 0.254 & \textbf{0.283} \\
    \hline
    Snare drum & \textbf{0.308} & 0.298 \\
    \hline
    Bass drum & 0.290 & \textbf{0.629} \\
    \hline
    \end{tabular}
    \vspace{0.2cm}
    \caption{
    \label{table:Eval-2-ENST}
    F-scores for ENST-Drums
    using the signal-adaptive median threshold evaluation.
    Higher F-scores are marked in bold.}
\end{minipage}
\hfill
\begin{minipage}{0.45\textwidth}
    \begin{tabular}{|c|c|c|}
    \hline
     & MUR & NeNMF \\
    \hline
    Snare drum & 0.275 & \textbf{0.369} \\
    \hline
    Bass drum & 0.444 & \textbf{0.454} \\
    \hline
    \end{tabular}
    \vspace{0.2cm}
    \caption{\label{table:Eval-2-song3}
    F-scores for ``Every Now and Then" using the signal-adaptive median threshold evaluation. Higher F-scores are marked in bold.}
\end{minipage}
\end{table}

\Cref{table:Eval-2-ENST,table:Eval-2-song3} present the F-scores evaluated using the signal-adaptive median threshold.
The scores are averaged over 28 tracks for ENST-Drums.
The F-scores are generally higher for NeNMF than for MUR.
The highest F-score of 62.9\% is achieved using NeNMF to detect the bass drum in ENST-Drums, and this F-score is over $2\times$ higher than the corresponding F-score of 29.0\% achieved using MUR.

\section{Conclusion}\label{sec:conclusion}
In our empirical results, NeNMF is consistently more accurate than MUR given a fixed computational budget. 
Further, NeNMF offers advantageous theoretical convergence guarantees. 
Therefore, we recommend using NeNMF for automatic drum transcription and other partially fixed NMF applications in the future.

As a limitation of this work, the evaluation procedure for the detected drum onsets is sensitive to parameters, and it is unclear which choice works best theoretically or empirically.
Further work on a practical and accurate evaluation method would enrich these experiments. 
As another limitation, this work focuses on two of the many possible optimizers for NMF.
Future work comparing other optimization methods or exploring different loss functions like the Kullback-Leibler or Itakura-Saito divergences could be interesting.

In conclusion, this work makes progress toward an accurate and interpretable automatic drum transcription system, first by highlighting partially fixed NMF as an interpretable model choice and second by identifying NeNMF as an efficient optimizer.
We acknowledge that the interpretable ADT system in this work does not yet achieve the ideal accuracy level, and indeed a super high-accuracy ADT system (e.g., with $99\%$ accuracy) is not yet available \cite{OST,musicanalysis,CNN,GillisBSS,SmaragdisConvNMF,Polyphonic,WuLerch,SemiSupCNMF,YooEtAl}.
Indeed, our validation efforts were limited by the amount of accessible drum audio data with accurate annotations, and we manually annotated over 600 drum hits for ``Every Now and Then''.
Yet, with further research on defining drum ``onsets'' and further modifications to the loss function, we believe it may be possible to build the partially fixed NMF approach into a fully accurate, interpretable ADT system so that manual annotations become unnecessary in the future.

\section*{Acknowledgments}
We would like to acknowledge Dan Manil, Zachary Lawrence, and Shawn Pana for collaboration in writing, performing, and recording ``Every Now and Then."


\appendix
\label{appendix}

\section{Proof of \texorpdfstring{\cref{thm:PFMUR}}{Thm.~3.2}}
\label{proof}
\begin{proof}
This section will prove that the update rule applied to $\bm{H}_D$ does not increase the loss function $F(\bm{H}_D, \bm{W}_H, \bm{H}_H)$.
A similar proof shows that the update rule applied to $\bm{W}_H$ or $\bm{H}_H$ does not increase the loss function.

We observe that the loss function is separable in the columns of $\bm{H}_D$:
\begin{equation*}
    F(\bm{H}_D, \bm{W}_H, \bm{H}_H)
    = \frac{1}{2} \sum_{j=1}^n \| \bm{V}(\cdot, j) - [\bm{W}_D\bm{H}_D(\cdot, j) + \bm{W}_H\bm{H}_H(\cdot, j)] \|_2^2.
\end{equation*}
We consider the change to the loss function after updating a single column from $\bm{H}_D(\cdot, j)$ to $\bm{x} \in \mathbb{R}^r$:
\begin{equation*}
    \tilde{F}(\bm{x}) = F(\bm{H}_D + [\bm{x} - \bm{H}_D(\cdot, j)] \bm{e}_j^\top, \bm{W}_H, \bm{H}_H).
\end{equation*}
$\tilde{F}(\bm{x})$ is a quadratic function of $\bm{x}$, and the gradient and Hessian terms are given as follows:
\begin{equation*}
    \nabla \tilde{F}(\bm{x})
    = \bm{W}_D^\top [\bm{W}_D\ \bm{x} + \bm{W}_H\bm{H}_H(\cdot, j) - \bm{V}(\cdot, j)], \qquad
    \nabla^2 \tilde{F}(\bm{x})
    = \bm{W}_D^\top \bm{W}_D.
\end{equation*}
Hence, $\tilde{F}(\bm{x})$ is specified by the formula
\begin{equation*}
\begin{aligned}
    \tilde{F}(\bm{x}) &= \tilde{F}(\bm{H}_D(\cdot, j)) + [\bm{x} - \bm{H}_D(\cdot, j)]^\top \nabla \tilde{F}(\bm{H}_D(\cdot, j)) \\
    & + \frac{1}{2} [\bm{x} - \bm{H}_D(\cdot, j)]^\top \bm{W}_D^\top \bm{W}_D [\bm{x} - \bm{H}_D(\cdot, j)].
\end{aligned}
\end{equation*}
To complete the proof, it suffices to show that updating the $j$th column from $\bm{H}_D(\cdot, j)$ to
\begin{equation}
\label{eq:update}
    \bm{x}_\star = \bm{H}_D(\cdot, j) \odot \frac{\bm{W}_D^\top\bm{V}(\cdot, j)}{\bm{W}_D^\top[\bm{W}_D\bm{H}_D(\cdot, j) + \bm{W}_H\bm{H}_H(\cdot, j)]}
\end{equation}
leads to a loss function $\tilde{F}(\bm{x}_\star) \leq \tilde{F}(\bm{H}_D(\cdot, j))$.

We introduce a diagonal matrix $\bm{K} \in \mathbb{R}^{r \times r}$ with elements
\begin{equation*}
    \bm{K}(i, i) = 
    \begin{cases}    
    \frac{\bm{W}_D(\cdot, i)^\top [\bm{W}_D \bm{H}_D(\cdot, j) + \bm{W}_H \bm{H}_H(\cdot, j)]}{\bm{H}_D(i, j)},
    & \bm{H}_D(i,j) > 0, \\
    0, & \bm{H}_D(i,j) = 0.
    \end{cases}
\end{equation*}
For any multiplicative update $\bm{x} = \bm{H}_D(\cdot, j) \odot \bm{y}$ with $\bm{y} \in \mathbb{R}^r$, we can upper bound
\begin{align*}
    & (\bm{H}_D(\cdot, j) \odot \bm{y})^\top \bm{W}_D^\top  \bm{W}_D (\bm{H}_D(\cdot, j) \odot \bm{y}) \\
    &= \sum_{i=1}^r \sum_{\ell = 1}^r \bm{y}(i) \bm{y}(\ell) \bm{H}_D(i, j) \bm{W}_D(\cdot, i)^\top \bm{W}_D(\cdot, \ell) \bm{H}_D(\ell, j) \\
    &\leq \sum_{i=1}^r \sum_{\ell = 1}^r \biggl[\frac{1}{2} \bm{y}(i)^2 + \frac{1}{2}\bm{y}(\ell)^2\biggr] \bm{H}_D(i, j) \bm{W}_D(\cdot, i)^\top \bm{W}_D(\cdot, \ell) \bm{H}_D(\ell, j) \\
    &= \sum_{i=1}^r \sum_{\ell = 1}^r \bm{y}(i)^2 \bm{H}_D(i, j) \bm{W}_D(\cdot, i)^\top \bm{W}_D(\cdot, \ell) \bm{H}_D(\ell, j) \\
    &\leq \sum_{i=1}^r \sum_{\ell = 1}^r \bm{y}(i)^2 \bm{H}_D(i, j) \bm{W}_D(\cdot, i)^\top [\bm{W}_D(\cdot, \ell) \bm{H}_D(\ell, j) + \bm{W}_H(\cdot, \ell) \bm{H}_H(\ell, j)] \\
    &= (\bm{H}_D(\cdot, j) \odot \bm{y})^\top \bm{K} (\bm{H}_D(\cdot, j) \odot \bm{y}).
\end{align*}
Therefore, for any multiplicative update $\bm{x} = \bm{H}_D(\cdot, j) \odot \bm{y}$, the quadratic function
\begin{equation*}
    G(\bm{x}) = \tilde{F}(\bm{H}_D(\cdot, j)) + [\bm{x} - \bm{H}_D(\cdot, j)]^\top \nabla \tilde{F}(\bm{H}_D(\cdot, j)) + \frac{1}{2} [\bm{x} - \bm{H}_D(\cdot, j)]^\top \bm{K} [\bm{x} - \bm{H}_D(\cdot, j)]
\end{equation*}
provides an upper bound on the loss function $\tilde{F}(\bm{x})$.

Last, the minimizer of $\bm{y} \mapsto G(\bm{H}_D(\cdot, j) \odot \bm{y})$ satisfies the gradient equal to zero condition:
\begin{equation*}
    0 = \bm{H}_D(\cdot, j) \odot \bm{K}[\bm{K}^+ \nabla \tilde{F}(\bm{H}_D(\cdot, j) + \bm{H}_D(\cdot, j) \odot \bm{y} - \bm{H}_D(\cdot, j)].
\end{equation*}
Equivalently, the minimizer satisfies the following identity:
\begin{align*}
    \bm{H}_D(\cdot, j) \odot \bm{y}
    &= \bm{H}_D(\cdot, j) - \bm{K}^+ \nabla \tilde{F}(\bm{H}_D(\cdot, j) \\
    &= \bm{H}_D(\cdot, j) - \frac{\bm{H}_D(\cdot, j) \odot  \bm{W}_D^\top [\bm{W}_D \bm{H}_D(\cdot, j) + \bm{W}_H\bm{H}_H(\cdot, j) - \bm{V}(\cdot, j)]}{\bm{W}_D^\top [\bm{W}_D \bm{H}_D(\cdot, j) + \bm{W}_H \bm{H}_H(\cdot, j)]} \\
    &= \bm{H}_D(\cdot, j) \odot \frac{\bm{W}_D^\top \bm{V}(\cdot, j)}{\bm{W}_D^\top [\bm{W}_D \bm{H}_D(\cdot, j) + \bm{W}_H \bm{H}_H(\cdot, j)]}.
\end{align*}
Hence, the minimizer of $G$ is precisely the multiplicative update rule $\bm{x}_\star = \bm{H}_D(\cdot, j) \odot \bm{y}$ given in \cref{eq:update}.
Since $G$ is an upper bound on $\tilde{F}$ which intersects at $\bm{H}_D(\cdot, j)$, we conclude:
\begin{equation*}
    \tilde{F}(\bm{x}_{\star}) \leq G(\bm{x}_{\star}) \leq G(\bm{H}_D(\cdot, j)) = \tilde{F}(\bm{H}_D(\cdot, j)).
\end{equation*}
This shows that the multiplicative update rule does not increase the loss function.
\end{proof}

\end{document}